\documentclass[9pt]{article}

\usepackage{amssymb, amsmath, latexsym, amscd, amsthm}
\usepackage{pdfsync}
\usepackage{graphicx, color, epsfig}  
\usepackage{epstopdf}
\usepackage[all,cmtip]{xy}
\usepackage{cite}

\def\sref#1{Section~\ref{#1}}

\def\tref#1{Theorem~\ref{#1}}
\def\dref#1{Definition~\ref{#1}}

\def\rref#1{Remark~\ref{#1}}
\def\cref#1{Corollary~\ref{#1}}
\def\pref#1{Proposition~\ref{#1}}
\def\lref#1{Lemma~\ref{#1}}
\def\eref#1{Example~\ref{#1}}

\newtheoremstyle{thm}{1em}{1em}{\itshape}{}{\bfseries}{.}{.5em}{} 
\newtheoremstyle{rem}{1em}{1em}{}{}{\bfseries}{.}{.5em}{}           

\theoremstyle{thm}                
\newtheorem{thm}{Theorem}[section]
\newtheorem{prop}[thm]{Proposition}

\newtheorem{cor}[thm]{Corollary}
\newtheorem{lem}[thm]{Lemma}

\newtheorem{dfn}[thm]{Definition}

\theoremstyle{rem} 
\newtheorem{rmk}[thm]{Remark}

\newtheorem{ex}[thm]{Example}

\newcommand{\intro}[1]
{\renewcommand{\thesection}{\fnsymbol{section}}
\setcounter{section}{-1}
\section{#1}
\renewcommand{\thesection}{\arabic{section}}
}

\DeclareMathOperator{\Ker}{Ker}
\DeclareMathOperator{\im}{Im}

\newcommand{\Pf}{\noindent{\bf Proof}. }
\newcommand{\cqfd}
{%
\mbox{}%
\nolinebreak%
\hfill%
\rule{2mm}{2mm}%
\medbreak%
\par%
}
 
\renewcommand{\b}{{\mathcal B}{}}

\newcommand{\CE}{{\mathcal C}{}}

\newcommand{\D}{{\mathcal D}{}}

\newcommand{\N}{{\mathbb N}{}}

\newcommand{\R}{\mathbb R}

\newcommand{\z}{{\mathcal Z}{}}
\newcommand{\Z}{\mathbb Z}

\newcommand{\ol}{\overline}
\newcommand{\ul}{\underline}

\newcommand{\rp}{respectively }

\begin{document} 

\title{A correspondence between contact structures and a class of cone structures}

\author{M\'elanie Bertelson\footnote{Chercheur Qualifi\'e F.N.R.S. ---
D\'epartement de Math\'ematique, Universit\'e Libre de Bruxelles --- Boulevard du Triomphe, 1050 Bruxelles,  Belgium --- {\tt mbertels@ulb.ac.be} --- Work supported by the Belgian Interuniversity Attraction Pole (IAP) within the framework ÒDynamics, Geometry and Statistical PhysicsÓ (DYGEST).}, C\'edric De Groote\footnote{Department of Mathematics, Stanford University, Stanford, CA 94305-2125, USA, {\tt cedricd@stanford.edu}}}

\date{\today}

\maketitle

\begin{abstract}

In his paper \emph{Cycles for the Dynamical Study of Foliated Manifolds and Complex Manifolds}, Denis Sullivan proves that a  closed manifold supports a symplectic structure if and only if it admits a distribution of cones of bivectors satisfying two conditions. We prove a similar result for contact structures. It relies on a suitable variant of the symplectization process that produces a $S^1$-invariant nondegenerate $2$-form on the closed manifold $S^1 \times M$ that is closed for a twisted differential. 
\end{abstract}

\intro{Introduction}

In \cite{S76}, Sullivan establishes, among other results, a correspondence, for closed manifolds, between symplectic structures and cone structures of bivectors satisfying certain hypotheses. A cone structure of $k$-vectors is a continuous distribution of compact convex cones $(C_x)_{x\in M}$, each $C_x$ contained in $\Lambda^kT_xM$. To such a cone structure is associated a cone ${\mathcal C}$ in the space of $k$-currents called the \emph{cone of structures currents}. It is defined to be the closed convex hull of the collection of all Dirac currents $\delta_P : \omega \mapsto \omega(P)$ associated to elements $P$ in the cone structure. When the manifold is closed, this cone is compact and plays an important role through the notion of positive differential forms, that is to say $k$-forms $\beta$ for which $c(\beta) > 0$ for each nontrivial structure current $c$. \\

A symplectic structure $\omega$ induces a cone structure through the choice of an almost complex structure $J$ tamed by $\omega$. It is ample and the associated cone of structure currents does not contain any nontrivial exact current (in the image of the adjoint of the exterior differential). Conversely, an ample cone structure $C$ on $M$ with no nontrivial exact structure current admits a contractible collection of positive differential forms that are symplectic. The proof uses the Hahn Banach separation theorem, as well as the duality between forms and currents. \\

The purpose of this paper, whose motivation is explained below, is to propose a correspondence between contact structures and certain cone structures. Our first idea was to apply Sullivan's correspondence to the symplectization of a contact manifold. This raises a number of difficulties. One is that the usual symplectization yields an open manifold~: $\R \times N$, while the assumption that $M$ is closed in Sullivan's correspondence is necessary because it implies that $\CE$ is compact, whence that Hahn Banach may be applied. Another issue is that, even if this obstacle could be circumvented, the symplectic structures we would recover on $\R \times N$ might not yield contact structures on $N$. Indeed, only $\R$-equivariant symplectic structures on $\R \times N$ naturally induce contact structures on $N$. These two problems can be bypassed by considering another version of the symplectization process that produces, from a contact form on $N$, a $S^1$-invariant nondegenerate $2$-form on $S^1 \times N$ that is closed for a twisted differential $D\beta = dt \wedge \beta + d \beta$. Such a form, called hereafter a \emph{$S^1$-invariant $D$-symplectic form}, yields, through the choice of a compatible $S^1$-invariant almost complex structure, an ample $S^1$-invariant cone structure that has no nontrivial $D^*$-exact structure current. Conversely, applying an invariant version of the Hahn Banach separation theorem to such a cone structure produces a $S^1$-invariant $D$-symplectic form on $S^1\times N$, whence a contact form on $N$. \\

It should specified that the cone structure associated to a certain cooriented contact structure depends on certain choices : that of a contact form and of an almost complex structure compatible with its symplectization. Nevertheless, the collection of induced cone structures is contractible. Likewise, given a ample cone structure of bivectors on $S^1 \times M$ that satisfies the homological condition, there is a large collection of associated contact forms (depending principally on the choice of closed hyperplane guaranteed by Hahn Banach) but, here too, the collection of possible choices is contractible. Now if the initial contact structure is not coorientable, the associated cone structure is supported by the manifold $S^1 \times \tilde{M}$, where $\tilde{M} \to M$ is the coorientation double cover of the contact structure, and is $\Z_2$-skew-invariant. \\

The twisted differential $D$ on $S^1\times M$ is in fact the differential of the Lichnerowicz cohomology associated to the closed form $dt$ and a $D$-symplectic form is a locally conformal symplectic (lcs) structure whose Lee form is $dt$. In view of these facts, it is tempting to believe that a similar result can be established for lcs structures and it is indeed the case. The correspondence is between lcs structures whose Lee form is a fixed closed form $\theta$ and ample cone structures of bivectors whose cones of structure currents do not contain non-trivial $D_\theta^*$-boundaries, where $D_\theta^*$ denotes the adjoint of the differential $D_\theta : \alpha \mapsto d \alpha + \theta \wedge \alpha$ for the Lichnerowicz cohomology associated to $\theta$. So the correspondence is really between locally conformal symplectic forms and ample cones of bivectors whose cones of structure currents avoid at least one of the various spaces of boundaries $\im D_\theta^*$. \\

Dusa McDuff has proven in her paper {\em Applications of convex integration ...} \cite{McD87} a criterion for existence of contact structures that also involves Sullivan's cone structures but, as explained hereafter, her approach is different. Starting from an orientable odd-dimensional manifold $M$ endowed with an exact $2$-form $d\alpha$ of maximal rank (or even a $2$-form of maximal rank, as the h-principle for odd-dimensional symplectic structures implies that these two data are equivalent up to homotopy), she searches for a closed $1$-form $\kappa$ such that $\alpha' = \alpha + \kappa$ is a contact form. It is the case if and only if the kernel of $d \alpha$ and that of $\alpha'$ are transverse, of course. Observe that, because $M$ is orientable, $\Ker d \alpha$ is a trivial line bundle and is therefore divided into two connected components by the $0$-section. McDuff considers the cone structure consisting of one of these components and formulates a condition on the associated cone of structure currents that is equivalent to existence of the desired $1$-form. This condition is that the collection of structure currents vanishing on $\alpha$ intersects trivially the (necessarily nontrivial) space of structure boundaries. It is also based on the Hahn Banach theorem. To summarize, McDuff starts from an exact $2$-form $d \alpha$ and builds a cone structure that depends on that $2$-form and which, when it satisfies a certain homological condition involving $\alpha$, implies existence of a contact structure differing from $\alpha$ by a closed $1$-form. So the cone structure and the homological condition it has to satisfy in order to induce a contact structure both depend on the initial choice of exact $2$-form. In particular she does not formulate general conditions on a cone structure that guarantee existence of a contact form. \\

Our motivation to transpose Sullivan's correspondence to the contact world originated in our interest for symplectic and contact structures of weaker than smooth regularity. There is a growing interest for the notion of $C^0$-symplectic structure, but other categories of symplectic and contact manifolds could also be considered, as, for instance, the PL or bi-Lipschitz ones. A fundamental problem is to compare these various categories. It is, for instance, tempting to believe that any PL symplectic or contact manifold of small dimension ($\leq 7$) may be smoothed. In order to prove that, one idea is to build from the PL symplectic or contact structure a cone structure of bivectors \`a la Sullivan satisfying the necessary conditions for inducing a smooth symplectic or contact structure. Now the reason for wanting a contact version of Sullivan's correspondence instead of testing this procedure in the symplectic category is simply that after the dimension $2$, for which most phenomena are quite simple, comes the dimension $3$, that belongs to the contact world. \\ 

The paper is organized as follows. The first section describes Sullivan's original construction. The second one presents a version of the symplectization process that is useful to us. The third section discusses the invariant Hahn Banach separation theorem. The fourth one gathers the results from the previous sections to establish the promised correspondence and the last one treats the non-coorientable case as well as locally conformal symplectic structures. 

\section*{Acknowledgments} We are extremely grateful to Dennis Sullivan for several very stimulating and inspiring conversations. We warmly thank Fran\c cois Laudenbach for helping us understand a delicate point about the topology of currents that we had overlooked. Finally, we are indebted to the anonymous referee for bringing to our attention a part of one of McDuff's papers related to our approach that we were unaware of.

\section{Sullivan's construction}

Let us recall some basic notions about convex sets.

\begin{dfn}\label{Felini} Let $V$ denote a topological vector space over the reals.
\begin{enumerate}
\item A cone in $V$ is a subset $\CE$ that is invariant under multiplication by positive reals. It is said to be convex when it is a convex subset of $V$.
\item A continuous linear form $\alpha \in V'$ is said to be positive on a cone $\CE \subset V$ if $\alpha(v) >0$ for each $v \in \CE-\{0\}$.
\item If a cone $\CE \subset V$ admits a positive linear form $\alpha$, the subset $\ul{\CE} = \alpha^{-1}(1) \subset \CE-\{0\}$ is called a base for $\CE$. The set  $\ul{\CE}$ is in bijective correspondence with the collection of rays in $\CE$. 
\item A cone is compact when it admits a positive linear form $\alpha$ and $\ul{\CE}$ is compact (notice that all bases are homeomorphic).
\end{enumerate}
\end{dfn}

The following definition is due to Sullivan's (cf.~\cite{S76}).

\begin{dfn}
A cone structure of $k$-vectors on a manifold $M$ is a continuous field $C = (C_x)_{x\in M}$ of compact convex cones $C_x \subset\Lambda^kT_xM$. 
\end{dfn}

To make sense of continuity, one chooses a Riemannian metric $g$ on $\Lambda^kTM$. It induces a distance $d$ on $S\Lambda^kTM$ the sphere bundle in $\Lambda^kTM$ and a Hausdorff distance $\rho$ on the collection $K$ of non-empty compact subsets of $S\Lambda^kTM$~:
$$\rho (S_1, S_2) = \max\Bigl\{\sup_{x \in S_1} d(x, S_2), \sup_{y \in S_2} d(S_1, y)\Bigr\}.$$
The cone structure $C$ can be seen as a map from $M$ to $K$ and continuity of $C$ can thus be defined by means of~$\rho$.

\begin{dfn} A $k$-form $\beta$ is said to be transverse to a cone structure $C$ of $k$-vectors, or positive on $C$, if $\beta_x(P) > 0$\footnote{If $P = v_1 \wedge ... \wedge v_k$, then $\beta(P) = \beta(v_1, ..., v_k)$.} for all $P \in C_x - \{0\}$ and all $x \in M$.
\end{dfn}

Due to compactness of the cones, transverse $k$-forms always exist. Indeed, for each $x$ in $M$, consider an element $\beta_x$ in the dual $\Lambda^kT^*_xM$ of $\Lambda^kT_xM$ such that $\alpha_x(P) >0$ for all non-vanishing $P \in C_x$. Because $C$ is continuous, the $k$-form $\beta_x$ can be extended to a positive $k$-form $\ol{\beta}_x$ on a neighborhood of $x$. Finally, the forms $\ol{\beta}_x$ may be glued into a global positive $k$-form by means of a partition of unity.

\begin{ex}\label{sympl-cone} Examples of cone structures that are essential here are the cone structures induced by symplectic forms. Let $\omega$ be such a form on a manifold $M$ and consider an almost complex structure $J$ compatible with the symplectic structure $\omega$\footnote{Which means that $\omega(Jv, Jw) = \omega(v,w)\; \forall \; v, w \in TM$ and $\omega(v, Jv) >0 \; \forall \; 0 \neq v \in TM$.}. Define the cone structure $C^J = {(C^J_x)}_{x\in M}$ whose fiber at $x$ is defined to be the convex closure of the cone 
$$\{v \wedge Jv \;|\; v \in TM\}.$$
The proof of \lref{CJ} shows that $C^J$ is indeed a cone structure. The original symplectic form $\omega$ is a positive form on $C^J$, as is any symplectic form that tames $J$. Alternatively, one could consider, instead of $C^J_x$, the convex closure of the collection of bivectors $v \wedge w \in \Lambda^2T_xM$ such that $\omega_x(v, w) > 0$, but that cone is not compact and its topological closure is too large in the sense that it contains isotropic bivectors as well. 
\end{ex}

Let us now recall some basic facts about currents that will be needed in the sequel. Let $M$ be a manifold. For $k \in \N$, a $k$-current on $M$ is a continuous linear form on the space $\Omega_c^k(M)$ of compactly supported $k$-forms on $M$ endowed with its standard Fr\'echet topology. The collection of $k$-currents, that is, the topological dual of  $\Omega_c^k(M)$, is denoted by $\D_k(M)$ and endowed with the strong topology. The transpose of the exterior differential $d : \Omega^k_c(M) \to \Omega^{k+1}_c(M)$ is a continuous differential $d^* : \D_{k+1}(M) \to \D_{k}(M)$.\\

It is a standard fact (see \cite{dR84} Th\'eor\`eme 13 p.~89) that $\Omega_c^k(M)$ is reflexive which means  that it is isomorphic to the strong dual of its strong dual through the evaluation map. In particular, any continuous linear form on currents corresponds to a unique differential form. Moreover, the forms vanishing on the space of exact (\rp closed) currents are the closed (\rp exact) forms. This is a consequence of reflexivity and of the following relation satisfied by the transpose of a continuous linear map $F$ between topological vector spaces (see \cite{T06} Formula (23.2) p~241)~:
\begin{equation}\label{ker-im}
\Ker F^* = (\im F)^\perp.
\end{equation}

Another key point is that the set of exact currents is a closed subspace of the space of all currents or equivalently the differential $d^*$ is a homomorphism\footnote{Let us recall that a continuous linear map $F : V \to W$ between topological vector spaces is said to be a homomorphism if the induced map $\ol{F} : V/\Ker F \to \im F$ is a homeomorphism for the obvious topologies.}. It is a consequence of the standard result recalled in \pref{TH} and of the closedness of the space of exact forms, itself a consequence of de Rham's theorem. The trivial observation that a subspace that is closed for the weak topology is also closed for the strong one is also needed.\\

Coming back to cone structures, to any such is associated a compact convex cone in the space of currents as described below. First observe that an element $P$ in some $\Lambda^kT_xM$ induces the Dirac current $\delta_P$ mapping a $k$-form $\beta$ to $\beta(P)$. Given a cone structure $C$ on $M$, consider the collection $\D_C = \{\delta_P \;|\;P \in C\}$ of Dirac currents. It is a (non-convex) cone in $\D_k(M)$.

\begin{dfn}
To a cone structure $C$ on $M$ is associated the cone of structure currents, that is, the topological closure of the convex closure of the cone $\D_C$. It is denoted by $\CE$. A closed current that belongs to $\CE$ is called a \emph{structure cycle}.
\end{dfn}
The subset $\CE$ is of course a closed convex cone and, when $M$ is a compact manifold, it is also compact (cf.~\cite{S76}, Proposition I.5.~p.~230). This fact seems more acceptable after having observed that a subset of $\D_k(M)$ is compact for the strong topology if and only if it is compact for the weak topology. 

\begin{rmk} We might as well consider the weak topology on the space of currents, as all the  properties we need are true for that topology as well.

\end{rmk}

Now one of Sullivan's main points in \cite{S76} is that certain properties of transverse $k$-forms may be encoded in the structure cone. For instance, whether exact or closed transverse forms exist is determined by the position of the structure cone relative to the subspaces of closed and exact currents, denoted respectively by $\z_k(M)$ and $\b_k(M)$. That relative position can be of either one of the following types~: 
\begin{enumerate}
\item[I] $\CE$ intersects $\b_k(M)$ non-trivially,
\item[II] $\CE$ intersects $\b_k(M)$ trivially and $\z_k(M)$ non-trivially,
\item[III] $\CE$ intersects $\z_k(M)$ trivially.
\end{enumerate}

Here is a precise statement~:

\begin{thm}\label{Sullivanthm}(Sullivan \cite{S76} Theorem I.7 p.~231)
Let $C$ denote a cone structure of $k$-vectors on a closed manifold $M$. Then 
\begin{enumerate}
\item There are either non-trivial closed structure currents or closed transverse forms.
\item If there are no closed transverse forms, then there are non-trivial exact structure currents (this corresponds to case I).
\item If there are no non-trivial closed structure currents, then there are exact transverse forms (this corresponds to case III).
\end{enumerate}
\end{thm}

More is said in \cite{S76} about case II, but we do not need this here. The main ingredient in the proof is the Hahn Banach separation theorem (see \cite{T06}, Proposition 18.2 p.~191 for instance), one version of which we recall below~: \\

\noindent
{\bf Hahn Banach separation theorem.} Let $W$ be a closed subspace in a  locally convex topological vector space $V$ that does not intersect a compact convex subset $C \subset V$. Then $W$ can be extended to a closed hyperplane $H$ that does not meet $C$ either.\\

The first point of \tref{Sullivanthm} is a consequence of the other two. The idea for the second point is that if the structure cone does not intersect $\b_k(M)$ (except along $0$), then, because $\CE$ is compact, the separating Hahn Banach theorem implies that the closed subspace $\b_k(M)$ is contained in a closed hyperplane $H$ that does not intersect a base for $\CE$. That closed hyperplane is the kernel of a continuous linear functional on $\D_k(M)$, determined up to a non-zero factor. By duality between forms and currents, that functional determines a $k$-form which can be chosen to be positive on $\CE$, whence on $C$ since $\D_C \subset \CE$, and which is closed because it vanishes on $\b_k(M)$. If $\CE$ is disjoint not only from $\b_k(M)$ but also from $\z_k(M)$, we obtain likewise exact transverse $k$-forms, proving the third point. \\

Coming back to the symplectic world, the question is how to encode the defining properties of a symplectic form into properties of an associated structure cone to which it is transverse. Recall the cone structure $C^J$ (introduced in \eref{sympl-cone}) associated to an almost complex structure $J$, itself compatible with a given symplectic structure $\omega$. The previous discussion shows that closedness of $\omega$ corresponds to absence of exact structure currents. How about nondegeneracy~? As explained hereafter, it corresponds to the property of ampleness for cone structures that is stated in \dref{ample} below. \\

As a preliminary, let us recall that the Schubert variety of a $2$-plane $\tau$ in a vector space $V$ is the collection, denoted by $S_\tau$, of all $2$-planes intersecting $\tau$ non trivially. Observe also that to a $2$-plane $\tau$ corresponds the line in $\Lambda^2V$ consisting of the bivectors $v \wedge w$ for $\{v, w\} \subset \tau$. This allows us to think of the Schubert variety $S_\tau$ as being a cone in $\Lambda^2V$.

\begin{dfn}\label{ample}
A cone $\CE$ in $\Lambda^2V$ is said to be ample if for any $2$-plane $\tau \subset V$, the Schubert variety of $\tau$ intersects $\CE$ non-trivially. A cone structure $C = {(C_x)}_{x\in M}$ is ample when each $C_x$ is ample in $\Lambda^2T_xM$.
\end{dfn}

\begin{lem}\label{cjample} The cone structure $C^J$ associated to an almost complex structure $J$ (cf.~\eref{sympl-cone}) is ample. Any $2$-form transverse to an ample cone structure of $2$-vectors is non-degenerate.
\end{lem}

\Pf For the proof that $C^J_x$ is ample, let $\tau$ be a $2$-plane in $T_xM$ and let $v \in \tau$. The bivector $v \wedge Jv$ belongs to $C^J_x$, which implies that $S_\tau$ intersects $C^J_x$ non-trivially. \\

For the second statement, let us suppose, on the contrary, that a $2$-form $\omega$ transverse to a cone structure $C$ is degenerate, that is, admits a $2$-plane $\tau$ in its radical. Let $\tau' \in S_\tau$ and let $\{v, w\}$ be a basis for $\tau'$ with $v \in \tau \cap \tau'$. We have $\omega(v, w) = 0$ which implies that the bivector $v \wedge w$ cannot belong to $C_x$, whence that $C_x \cap S_\tau = \emptyset$. Thus $C_x$ is not ample.

\cqfd

Now \tref{Sullivanthm} and \lref{cjample} imply the following proposition.

\begin{prop} (Sullivan, \cite{S76} Theorem III.2, p.~249) Let $M$ be a closed manifold. If a symplectic structure $\omega$ is given on $M$, then any choice of compatible almost complex structure $J$ induce an ample cone structure $C^J$ that has no exact structure cycles and admits $\omega$ as a transverse form. Conversely, an ample cone structure $C$ without exact structure cycles admits a non-empty contractile collection of transverse forms that are symplectic. 
\end{prop}

The statement about contractibility is easily seen. Indeed, the collection of forms that are positive on the cone structure and vanish on $\b^k(M)$ is convex.

\section{An invariant version of the symplectization process}\label{IVSP}

The symplectization of a contact manifold $(M, \xi)$ consists of the open manifold $\R\times M$ endowed with an $\R$-equivariant symplectic form. We present, in this section, a variant of the symplectization process which yields a bijective correspondence between contact forms on a closed manifold $M$ and $S^1$-invariant non degenerate $2$-forms on the manifold $S^1 \times M$ that are closed for a twisted differential.\\

Let us first recall the usual symplectization process. If  $(M, \alpha)$ is a closed cooriented contact manifold then $(\R \times M, d(e^s \pi^*\alpha))$, where $\pi$ is the canonical projection of $\R \times M$ onto $M$ and $s$ is the function $\R \times M \to \R : (s, x) \mapsto s$, is an open symplectic manifold whose symplectic form is equivariant with respect to the standard action $\rho (t, (s, x)) = \rho_t(s, x ) = (t + s, x)$ of $\R$ on $\R\times M$. That is, the form $\omega = d(e^s \pi^*\alpha)$ satisfies the relation
$$\rho_t^*\omega = e^t \omega, \;\forall t \in \R.$$

Conversely, an equivariant symplectic form on $\R\times M$ yields a contact form. Indeed, an equivariant $2$-form on $\R\times M$ is of the type
$$\beta = e^s \pi^*\beta_0 + e^s ds\wedge \pi^*\alpha_0,$$
for a $2$-form $\beta_0$ and $1$-form $\alpha_0$ on $M$. It is closed if and only if $d \alpha_0 = \beta_0$ and non-degenerate exactly when $\alpha_0 \wedge \beta_0 \wedge ... \wedge \beta_0$ does not vanish. \\

This yields a one-to-one correspondence between the collection ${\rm Cont}(M)$ of contact forms on $M$ and the collection ${\rm ESymp}(\R \times M)$ of equivariant symplectic forms on $\R \times M$. Since the aim is a correspondence between contact forms and forms on $S^1 \times M$,  we are interested in invariant rather than equivariant forms on $\R\times M$, as they pass to the quotient. 

\begin{rmk}
An invariant $2$-form on $\R\times M$, that is, a form $\beta \in \Omega^2(\R \times M)$ that satisfies $\rho_t^*\beta = \beta$ for all $t \in \R$, induces a $2$-form $\beta_0$ and $1$-form $\alpha_0$ on $M$ such that
$$\beta = \pi^*\beta_0 + ds\wedge \pi^*\alpha_0.$$
Notice that $\beta$ is closed if and only if both $\beta_0$ and $\alpha_0$ are closed. 
\end{rmk}

Thus an invariant symplectic $2$-form is far from inducing a contact form. On the other hand, an equivariant symplectic form induces a non-degenerate $2$-form  that is closed for a twisted differential. Indeed, observe that the map $\varphi : \beta \mapsto e^s \beta$ induces a one-to-one correspondence between invariant and equivariant forms on $\R \times M$. The property of non-degeneracy is of course preserved. Moreover, 
$$d (e^s \beta) = e^s (ds \wedge \beta + d\beta).$$
This means that the map $\varphi$ intertwines the exterior differential $d$ with the differential operator 
$$D : \Omega^{i}(\R\times M) \to \Omega^{i+1}(\R \times M) : \beta \mapsto D\beta = ds\wedge \beta + d\beta.$$
With symbols~:
$$d \circ \varphi = \varphi \circ D.$$

\begin{lem}\label{Dprop}
The operator $D$ obviously satisfies the following properties~:

\begin{enumerate}
\item[-] $D^2 = 0$,
\item[-] $\rho_t^* \circ D = D \circ \rho_t^*$, for all $t \in \R$,
\item[-] $D$ is continuous.
\end{enumerate}
\end{lem}

As mentioned in the introduction, the operator $D$ is the differential for the Lichnerowicz cohomology associated to the closed form $ds$. Since it commutes with the $\R$-action, it passes to the quotient under the $\Z$-action and yields a  differential operator, also denoted by $D$, on $\Omega^*(S^1 \times M)$. As a side remark, the latter is not a derivation of the wedge product. Instead, the following relation holds true~:
$$D(\alpha \wedge \beta) = D\alpha \wedge \beta + (-1)^{|\alpha|} \alpha \wedge d\beta.$$

Introducing the collection ${\rm ISymp}^D(\R \times M)$ of non-degenerate invariant $2$-forms $\beta$ on $\R\times M$ satisfying $D\beta = 0$, called hereafter $\R$-invariant $D$-symplectic forms, the previous discussion implies that the following map is a bijection~: 
$$\psi : {\rm Cont}(M) \to {\rm ISymp}^D(\R \times M) : \alpha \mapsto e^{-s} d(e^s \pi^*\alpha) = ds \wedge \pi^*\alpha + d(\pi^*\alpha).$$

Now we would like to pass from $\R \times M$ to $S^1 \times M$. The map $p : \R\times M \to S^1 \times M : (s, x) \mapsto (e^{is}, x)$ induces a push-forward $p_*$ from the space of $\R$-invariant forms on $\R \times M$ to the space of $S^1$-invariant forms on $S^1\times M$ defined by $$[p_*(\beta)]_{(e^{is}, x)} (v_1, ..., v_\star) = \beta_{(s,x)}(\ol{v}_1, ..., \ol{v}_\star),$$
where $\ol{v_i} \in T(\R \times M)$ is the lift of $v_i$ through the point $(s, x)$. It is of course independent on the choice of $(s, x)$ in $p^{-1}(e^{is}, x)$. \\

We have thus obtained a bijective correspondence between the collection of contact forms on $M$ and the collection, denoted by ${\rm ISymp}^D(S^1\times M)$, of $S^1$-invariant $D$-symplectic forms on $S^1 \times M$. More precisely, the map 

$${\mathcal S} : {\rm Cont}(M) \to {\rm ISymp}^D(S^1\times M) : \alpha \mapsto {\mathcal S}(\alpha) = p_*\bigl(e^{-s} d(e^s \pi^*\alpha) \bigr)$$
is one-to-one.

\begin{rmk}\label{closedexact} The ``twisted symplectization" ${\mathcal S}(\alpha)$ of a contact form $\alpha$ coincides with $D\alpha$. It is thus a $D$-exact form. It is easy to verify directly that a $D$-closed nondegenerate $S^1$-invariant $2$-form $\beta$ on $S^1 \times M$ is necessarily the $D$-boundary of a contact form on $M$. Indeed, since $\beta$ is invariant, it is of the type $\beta = \pi^*\beta_0 + dt \wedge \pi^*\alpha_0$, for some forms $\beta_0$ and $\alpha_0$ on $M$. The relation $D\beta = 0$ implies $\beta_0 = d \alpha_0$. Thus $\beta = D\pi^*\alpha_0$. Now $\beta^n = n dt \wedge \pi^*(\alpha_0 \wedge d\alpha_0^{n-1})$ which implies that $\alpha_0$ is a contact form if and only if $\beta$ is non degenerate.
\end{rmk}

\begin{rmk} The usual symplectization of a coorientable contact structure $\xi$ can be defined intrinsically as the subbundle of the cotangent bundle of $M$ consisting of the $1$-forms whose kernel contains the contact structure, endowed with the restriction of the Liouville symplectic form. Non vanishing sections of that subbundle and contact form for $\xi$ are one and the same thing. In comparison, the form ${\mathcal S}(\alpha)$ is not intrinsic. Indeed, while a contact form $\alpha$ is not intrinsic, the forms $e^s \alpha$ and therefore $d(e^s \alpha)$ on $\R \times M$ are intrinsic but the form $e^{-s} d(e^s \alpha)$ is not anymore. This observation does not affect much our point. It is true that it makes the cone structure associated to a contact structure (see \sref{CCFCD}) depend on the choice of a contact form but the cone structure depends anyway on the choice of an almost complex structure compatible with the symplectization. Moreoever, the collection of all possible cone structures induced by a fixed contact structure is contractible.
\end{rmk}

\section{An invariant version of the Hahn Banach separation theorem}

In this section appears a proof of the invariant separation Hahn Banach theorem that relies on the invariant analytic Hahn-Banach theorem. Surprisingly, a proof of the precise statement we need is not so easily accessible in the literature whence our decision to include it in the text. See nevertheless \cite{L74} Theorem~1 in Section~3, as well as \cite{DG14} Th\'eor\`eme 2.20 p~32.\\

Let us first recall the statement, due to R.~Agnew and A.~Morse \cite{A-M38}, of the invariant analytic Hahn Banach theorem. 

\begin{thm} (Agnew \& Morse) Let $V$ denote a real topological vector space and $W$ a linear subspace. Let $p : V \to \R^+$ be a positively homogeneous subadditive functional\footnote{This means that $p(\lambda v) = \lambda p(v) \; \forall \lambda \in \R^+, v \in V$ and $p(v+w) \leq p(v) + p(w) \; \forall \; v, w \in V$.} and let $f$ be a linear form on $W$ such that $f(w) \leq p(w)$ for all $w \in W$. If $G$ is a solvable group acting continuously on $V$, preserving $W$ and leaving $p$ and $f$ invariant, then there exists a linear extension of $f$ to $V$ that is invariant under $G$ and satisfies $f(v) \leq p(v)$ for all $v \in V$.
\end{thm}

Here is the version of the invariant geometric Hahn-Banach theorem used hereafter~:
 
\begin{thm}\label{invHB} Let $G$ be a solvable group acting continuously on a real Hausdorff topological vector space $V$. Suppose $K$ is a compact convex subset of $V$ invariant under $G$ and $W$ is a closed linear subspace of $V$ also invariant under $G$ and that does not meet $K$. If $K$ contains a fixed point $c_0$ for the action of $G$ and if the complement of $W - K$ contains a convex invariant open neighborhood of $0$, then the subspace $W$ can be extended to a closed invariant hyperplane that does not meet $K$.
\end{thm}

\Pf Let $A$ be a convex invariant open neighborhood of $0$ contained in $W - K$. Consider the subset $O = A+K-c_0$. It is convex, invariant and open as a union $\displaystyle{O = \bigcup_{c \in K} (A + (c-c_0))}$ of open sets. Consider its gauge 
$$p_O : V \to \R^+ : v \mapsto \inf\{\lambda >0 \;|\; v \in \lambda O\}.$$ 
The functional $p_O$ is positively homogeneous subadditive and invariant because $O$ is invariant. Moreover $O = p_O^{-1}([0,1))$ since $O$ is open. Thus $p_O$ is continuous. Now consider the subspace $W'$ generated by $W$ and $c_0$ and the linear functional $\phi : W' \to \R : w -t c_0 \mapsto t$. Observe that $W'$ is invariant and that $\phi$ is continuous (because $W'$ in closed) and invariant. Moreover, the relation $\phi(w') \leq p_O(w')$ holds true for all $w' \in W'$. Indeed, for $t>0$, we have
$$p_O(w - t c_0) = \displaystyle{t p_O(\frac{w}{t} - c_0)} \geq t$$
since $\displaystyle{\frac{w}{t} - c_0 \in W-K}$ which is disjoint from $O$.\\

The Agnew-Morse theorem above implies that $\phi$ can be extended to an invariant linear functional $\ol{\phi} : V \to \R$ such that $\ol{\phi}(v) \leq p_O(v)$ for all $v \in V$. Observe that, as implied by the previous inequality and the continuity of $p_O$, the functional $\ol{\phi}$ is continuous. Let $H$ denote the kernel of $\ol{\phi}$. It is a closed invariant hyperplane containing $W$. Moreover $H \cap K = \emptyset$ because if $c \in K$, then $c - c_0 \in O$ and thus $p_O(c-c_0) <1$. This implies that $$\ol{\phi}(c) = \ol{\phi}(c_0) + \ol{\phi}(c -c_0) \leq \phi(c_0) + p_O(c-c_0) < -1 + 1 = 0,$$
showing that $c$ cannot belong to $H$.
\cqfd

Under the additional assumptions that $G$ is compact and $V$ is locally convex, the awkward assumption of existence of a convex invariant open neighborhood of $0$ contained in $W - K$ is satisfied. To prove this we need the following standard result about actions of topological groups.

\begin{lem}\label{theta}
Consider a continuous action of a topological group $G$ on a topological space $X$. Then, for a compact subspace $\Theta$ of $G$ and a closed subset $C$ of $X$, the set $\Theta \cdot C = \{\theta \cdot x \; | \; \theta \in \Theta, x \in C\}$ is closed. 
\end{lem} 

\begin{cor}\label{cor} Let $G$ be a compact solvable topological group acting continuously on a real locally convex Hausdorff topological vector space $V$. Suppose $K$ is a compact convex subset of $V$ invariant under $G$ and $W$ is a closed linear subspace of $V$ also invariant under $G$ and that does not meet $K$. If $K$ contains a fixed point for the action of $G$ then the subspace $W$ can be extended to a closed invariant hyperplane that does not meet $K$.
\end{cor}

\Pf Observe that, under the assumptions of \cref{cor} the existence of the convex invariant open neighborhood of $0$ in the complement of $W - K$, needed in \tref{invHB}, is guaranteed. Indeed, the subset $W - K$ is closed, as implied by \lref{theta} applied to the additive action of $V$ on itself. Its complement contains therefore a convex open neighborhood $A_0$ of $0$. Then \lref{theta}, again, implies that the invariant subset $A = \cap_{g \in G} (g \cdot A_0)$ is open (since its complement $G \cdot (V - A_0)$ is closed). It is of course convex as well. 
\cqfd

\section{Correspondence contact forms/cone structures}\label{CCFCD}

We are now ready to establish a correspondence between contact forms and a certain type of cone structures. Recall from \sref{IVSP} that we may think of contact forms on $M$ as being $S^1$-invariant $D$-symplectic forms on $S^1 \times M$. The idea now is to extend Sullivan's result to forms and cone structures that are invariant under an $S^1$-action. \\

Set $P = S^1 \times M$. For an element $\omega$ in ${\rm Symp}^{S^1}_D(P)$, consider an invariant almost complex structure $J$ that is compatible with $\omega$. There is a natural class of such structures. Indeed, consider the contact form $\alpha$ associated with $\omega$ and its Reeb vector field $R_\alpha$ (defined by the relations $d\alpha(R_\alpha, \cdot) \equiv 0$ and $\alpha(R_\alpha) \equiv 1$). Choose a compatible almost complex structure $J_0$ on the symplectic vector bundle $(\xi = \Ker \alpha, d \alpha|_{\xi\oplus\xi})$. Now define a vector bundle morphism $J : TP \to TP$ by setting 
$$\left\{\begin{array}{lll}
\displaystyle{\frac{}{}J(v)} & = & J_0(v) \;\; \mbox{for} \;\;v \in \xi \\
J(R_\alpha) & = & \displaystyle{\frac{\partial}{\partial t}}.
\end{array}\right.$$
The morphism $J$ is an $S^1$-invariant almost complex structure compatible with $\omega$, as is easily verified.\\

Now we construct the field $C^J$ of cones of $2$-vectors on $P$ associated to the almost complex structure $J$ as is done in \rref{sympl-cone}. The cone at $p$ is thus the convex hull of the set 
$$\{v\wedge Jv \;|\; v \in T_pP \}.$$ 

\begin{lem}\label{CJ} The field $C^J$ is continuous, compact, ample and invariant.
\end{lem}
The proof is elementary but is nevertheless included. \\

\Pf The fact that $C^J$ is continuous follows directly from the fact that one can find local trivializations of the vector bundle $TP$ relative to which the almost complex structure $J$ and thus also the cone field $C^J$ are constant. \\

To prove compactness of $C$, let $p \in P$ and consider pointwise linearly independent local sections $e_1, ..., e_{2n}$ of $TP$ near $p$ such that $J(e_i) = e_{i+n}$, $i = 1, ...,n$. Define the local $2$-form
$$\beta = \sum_{i=1}^n e_*^i \wedge e_*^{i+n},$$
where $\{e_*^1, ..., e_*^{2n}\}$ is the basis of $T^*P$ dual to $\{e_1, ..., e_{2n}\}$. Then
$$\beta(v \wedge Jv) = \sum_{i=1}^{2n} v_i^2.$$ 
It is now obvious that $\beta_p^{-1}(1) \cap C_p\subset \Lambda^2T_pP$ is compact. \\

Ampleness of $C^J$ has been verified in the proof of \lref{cjample} and invariance follows directly from that of $J$.
\cqfd

\begin{lem} When $M$ is compact, the structure cone $\CE^J$ associated to $C^J$ is compact and invariant. 
\end{lem}

\Pf For the compactness of $\CE^J$, we refer to \cite{S76}. Invariance follows directly from that of $C^J$. 
\cqfd

It is now necessary to prove that the space of currents that are exact for the adjoint $D^*$ of the operator $D$, introduced in \sref{IVSP}, is a closed subspace of $\D_*(M)$. \\

Let us denote by $\z^k_D(P)$ (\rp $\b^k_D(P)$) the subspace of $\Omega^k(P)$ consisting of $D$-closed (\rp $D$-exact) $k$-forms and by $\z_k^{D^*} (P)$ (\rp $\b_k^{D^*}(P)$) the space of closed (\rp exact) currents for $D^*$. Because $\rho_t^* \circ D = D \circ \rho_t^*$ for all $t \in S^1$ (cf.~\lref{Dprop}), all those spaces  are $S^1$-invariant. 

\begin{lem}\label{bclosed} The space $\b^{D^*}_k(P)$ is a closed subspace of $\D_k(P)$. 
\end{lem}

\Pf First show that $\b^{k+1}_D(P)  = \im D$ is a closed subspace of $\Omega^{k+1}(P)$ and then invoke the following classical result (e.g.\cite{T06} Proposition 35.7 p. 366)~:

\begin{prop}\label{TH} For a continuous linear map $u$ between two real locally convex Hausdorff topological vector spaces $E$ and $F$, the following properties are equivalent 
\begin{enumerate}
\item[-] $u(E)$ is closed in $F$;
\item[-] the transpose $u^*$ of $u$ is a homomorphism of $F'$ onto $u^*(F') \subset E'$ when $F'$ and $E'$ are endowed with their weak topology.
\end{enumerate}
\end{prop}

Supposing that $\b^{k+1}_D(P)$ is closed and applying this proposition to $D : \Omega^k(P) \to \Omega^{k+1}(P)$ implies that $D^*$ is a homomorphism. Therefore, because the space $\D^k(P)/\Ker D^*$ is complete, the space $\im D^*$ is also complete, whence closed. This is true for the weak topology and thus for the strong topology too. \\

It remains now to show that $\b^k_D(P)$ is closed for all $k$. Recall from \sref{IVSP} the map
$$\varphi : \Omega^k(\R\times M) \to \Omega^k(\R \times M) : \beta \mapsto e^s \beta.$$
It is a homeomorphism that intertwines $d$ with $D$ and, therefore, that induces a bijection between their respective images $\b^k_D(\R \times M)$ and $\b^k(\R \times M)$. Since the space of exact forms is closed in the space of differential forms (even if the underlying manifold is not closed)
, this implies that $\b^k_D(\R \times M)$ is closed in $\Omega^k(\R \times M)$. Besides, the map 
$$p^* : \Omega^k(S^1 \times M) \to \Omega^{k}(\R \times M) : \beta \mapsto p^*(\beta).$$
induces a homeomorphism between $\Omega^k(S^1 \times M)$ and $\Omega^k(\R \times M)^\Z$, the set of fixed elements of the action of $\Z$ on $\Omega^*(\R \times M)$. The latter being continuous, the set $\Omega^k(\R \times M)^\Z$ is a closed subspace of $\Omega^k(\R \times M)$. Thus $\b^*_D(P)$ appears to correspond, under $p^*$, to an intersection $\Omega^k(\R \times M)^\Z \cap \b^*_D(\R \times M)$ of closed subspaces of $\Omega^k(\R \times M)$.     
\cqfd

Now to apply the invariant version of the Hahn Banach separation theorem, we need to have a fixed point for the action of $S^1$ on a basis for the structure cone. Let us recall Tychonoff's fixed point theorem.

\begin{thm} Let $E$ be a locally convex topological vector space, let $C$ be a compact convex subset of $E$ and let $f : C \to C$ be a continuous map. Then $f$ has a fixed point.
\end{thm}

Consider and element $g \in S^1$ that generates a dense subgroup in $S^1$. Tychonoff's fixed point theorem implies that $\rho_g : \ul{\mathcal C} \to \ul{\mathcal C}$ has a fixed point $c$. Now $c$ is fixed under the action of the subgroup generated by $g$ as well. Since that subgroup in dense in $S^1$ and the action of $S^1$ on $\ul{\mathcal C}$ is continuous, all elements of $S^1$ fix $c$. \\

Alternatively, one may use the Markov-Kakutani fixed point theorem whose statement is recalled hereafter.

\begin{thm} (Markov \cite{M36} and Kakutani \cite{K38}) Let $C$ be a compact convex subset of a Hausdorff topological vector space $E$ and let $G$ be a collection of  commuting continuous affine transformations of $E$ that preserve $C$. Then there exists a point in $C$ that is fixed under all elements of $G$.
\end{thm}

Existence of a fixed point allows us to apply the geometric Hahn-Banach theorem. 

\begin{thm}\label{thm} Let $M$ be a closed manifold. A cooriented contact structure $\xi$ on $M$ induces a non-empty contractible collection of ample $S^1$-invariant cone structures on $P = S^1 \times M$ with no non-trivial $D^*$- exact structure currents. Conversely, an ample $S^1$-invariant cone structure on $P$ with no non-vanishing $D^*$- exact structure current induces a non-empty contractible collection of contact structures on $M$.
\end{thm}

\Pf Given a cooriented contact structures $\xi$ there exists a contractible collection of forms $\alpha$ defining $\xi$ and compatible with its coorientation. Consider an invariant almost complex structure $J$ on $P$ compatible with the symplectization $\omega = {\mathcal S}(\alpha)$ of some contact form $\alpha$ defining $\xi$ (cf.~\sref{IVSP}) and the associated cone structure $C^J$. It is ample, invariant (cf.~\lref{CJ}) and the associated cone structure $\CE^J$ does not contain any $D^*$- exact structure cycle. Indeed, the form $\omega$ is positive on $\CE^J$ but a $D$-closed form may not be positive on $D^*$- exact currents. \\

Conversely, let $C$ denote an ample $S^1$-invariant cone structure on $P$ without non-vanishing $D^*$- exact structure current. Consider an invariant basis $\underline{\CE}$ for the structure cone $\CE$. To show that such a base exists, it suffices to construct an invariant positive form, which is easily done as follows. Let $\beta$ be a postitive form on $C$. Define $\beta'$ to be the invariant extension of the restriction of $\beta$ to $T_{\{0\} \times M}P$. More explicitly~:
$$\beta'_{(t,x)} = \rho_{-t}^*\bigl(\beta_{(0,x)}\bigr).$$
The form $\beta'$ remain positive because $C$ is invariant. Now Tychonoff's fixed point theorem implies existence of a fixed point $c$ in $\underline{\CE}$ for the action of $S^1$ and thus \cref{cor} implies that the closed subspace $\b^{D^*}_2(P)$ may be extended to a closed hyperplane that does not meet $\underline{\CE}$. That closed hyperplane is the kernel of a continuous linear functional $\alpha$ on $\D_2(P)$ positive on $\underline{\CE}$. The presence of the fixed point $c$ implies that $\alpha$ is invariant. The space $\Omega^2(P)$ being reflexive, that functional is induced by a $2$-form $\omega$ which is invariant and $D$-closed because it vanishes on $\b_2^{D^*}(P)$ (argument identical to the one that shows that the closed forms are the ones vanishing on exact currents and that uses reflexivity of $\Omega^k(P)$ together with formula (\ref{ker-im})). Finally, such a form is the symplectization ${\mathcal S}(\alpha)$ of a contact form $\alpha$ on $M$. 
\cqfd

\section{The non coorientable case}
 
If $\xi_0$ is a non coorientable contact structure on a manifold $N$, consider the coorientation double cover $p : M \to N$, its $\Z_2$-action $\rho : \Z_2 \times M \to M$ given by $\rho_a(x) = -x$ and the associated $\Z_2$-invariant coorientable contact structure $\xi = \pi_*^{-1}\xi_0$. The latter admits a $\Z_2$-skew-invariant defining $1$-form $\alpha$\footnote{The expression \emph{an object ${\mathcal O}$ is $\Z_2$-skew-invariant} means that, if $a$ denotes the nontrivial element of $\Z_2$, the action of $a$ on the object ${\mathcal O}$ produces $-{\mathcal O}$}. The space $P = S^1\times M$ inherits a $\Z_2$-action as well and, as explained hereafter, it supports a $Z_2$-skew-invariant $S^1$-invariant almost complex structure compatible with the twisted symplectization $D\alpha$ of $\alpha$. \\

We first claim that there exists an almost complex structure $J$ on $\xi$ which is compatible with $d\alpha$ and such that $\rho_a \cdot J = -J$, where, of course, $\rho_a \cdot J = (\rho_a)_* \circ J \circ (\rho_a)_*$. Indeed, one way to construct a compatible almost complex structures $J$ on $\xi$ is as follows~: choose a Euclidean structure $g$ on $\xi$, consider the bundle morphism $A : \xi \to \xi$ defined by $g(X, Y) = d\alpha(AX, Y)$, observe that it is skew-adjoint ($A^* = -A$), whence that $-A^2 = A^*A$ is symmetric and positive definite. Then define the almost complex structures $J = AQ^{-1}$, where $Q = \sqrt{-A^2}$. Now if $g$ is $\Z_2$-invariant, which is easily achieved since Euclidean structures may be averaged, the morphism $A$ is $\Z_2$-skew-invariant and therefore the almost complex structures $J$ too~:
$$\rho_a \cdot J = -J.$$
observe that we could not hope for a $\Z_2$-invariant almost complex structure compatible with $d\alpha$ here as it would yield an almost complex structure on $\xi_0$, which is impossible since the bundle $\xi_0$ is not a symplectic vector bundle in the non-coorientable case.\\

Now we would like to extend $J$ to the entire twisted symplectization in such a way that it remains $\Z_2$-skew-invariant. Observe that the Reeb vector field $R_\alpha$ is $\Z_2$-skew-invariant~: $(\rho_a)_{*} (R_\alpha) = - R_\alpha$. Now in order to obtain a $\Z_2$-skew-invariant almost complex structure $J$ on $P$, compatible with $D\alpha$, and such that $J(R_\alpha) = \pm \partial_s$, it must be that $(\rho_a)_* \circ J (R_\alpha) = J (R_\alpha)$. Indeed,
$$-J(R_\alpha) = (\rho_a \cdot J) (R_\alpha) = (\rho_a)_{*} \circ J (-R_\alpha) = (\rho_a)_{*} (-J(R_\alpha)).$$ 
So if we lift the $\Z_2$-action to $P = S^1 \times M$ as follows~:
$$\rho_a : P \to P : (z, x) \mapsto (z, -x),$$
the almost complex structure $J$ is $\Z_2$-skew-invariant.\\

Now the cone structure $C^J$ on $S^1 \times M$ associated to $J$ is ample, $S^1$-invariant and $\Z_2$-skew-invariant since
$$\begin{array}{lll}
\rho_a \cdot (v \wedge Jv) & = & (\rho_a)_* (v) \wedge (\rho_a)_* (Jv) \\
& = & (\rho_a)_* (v) \wedge (\rho_a)_* \circ J \circ (\rho_a)_* \circ (\rho_a)_* (v) \\
& = & - (\rho_a)_* (v) \wedge J \circ (\rho_a)_* (v).
\end{array}$$
Notice that we are not saying that the restriction of $\rho_a$ to $C^J$ coincides with the map $-Id$. Likewise the associated cone ${\mathcal C}^J$ of structure cycles is $S^1$-invariant and $\Z_2$-skew-invariant. Moreover it does not contain $D^*$-exact structure cycles. \\

Conversely, given an ample $S^1$-invariant, $\Z_2$-skew-invariant ample cone structure $C$ on $P$ with no $D^*$-exact structure cycles, one recovers $S^1$-invariant, $\Z_2$-skew-invariant $D$-symplectic forms on $P$, themselves inducing non-coorientable contact structures on $N$. The first step is to construct a $S^1$-invariant compact base $\ul{\mathcal C}$ for ${\mathcal C}$ that is also $\Z_2$-skew-invariant. Given an $S^1$-invariant form $\beta_0$ that is positive on $C^J$, the form 
$$\beta = \beta_0 - (\rho_a)^*\beta_0$$ 
is also positive on $C^J$ and thus $\beta^{-1}(1)$ yields the desired basis. It is now useful to consider the $\Z_2$-action on $\D_*(M)$ defined by~:

$$\eta_a : \D_*(M) \to \D_*(M) : c \to -\rho_a \cdot c.$$
This action commutes with that of $S^1$ and leaves both $\ul{\mathcal C}$ and $\b_2^{D^*}(P)$ invariant. The Markov-Kakutani fixed point theorem implies that $\ul{\mathcal C}$ contains a fixed point for the $\Z_2 \times S^1$-action. So the invariant Hahn Banach theorem, or rather \cref{cor}, implies that the space of $D^*$- boundaries is contained in a closed $\Z_2 \times S^1$-invariant hyperplane $H$ disjoint from $\ul{\mathcal C}$. Let $\omega$ denote a linear form on $\D_2(P)$ whose kernel is $H$ and that is positive on $\ul{\mathcal C}$. The form $\omega$ is in fact a $D$-closed differential form that is $\Z_2 \times S^1$-invariant and non-degenerate. \\

Now because $\omega$ is $D$-closed and $S^1$-invariant, it is the $D$-boundary of a contact form $\alpha$ (cf.~\rref{closedexact}). The fact that $\omega$ is $\eta$-invariant, that is $\rho_a^*\omega = - \omega$, implies the corresponding relation for $\alpha$~:
$$\rho_a^*\alpha = - \alpha.$$
Thus $\alpha$ yields a non coorientable contact structure on $N$. We have thus proven the following result.

\begin{thm} Let $\xi$ be a non-coorientable contact structure on a closed manifold $N$. Consider the coorientation double cover $\pi : M \to N$ associated to $\xi$ and the $\Z_2 \times S^1$-action on $P = S^1 \times M$ given by $\rho_a \cdot (s, x) = (s, -x)$. Then $\xi$ induces a non-empty contractible collection of ample $S^1$-invariant, $\Z_2$-skew-invariant cone structures on $P$ with no non-trivial $D^*$- exact structure currents. Conversely, an ample $S^1$-invariant, $\Z_2$-skew-invariant cone structure on some $P= S^1 \times M$ for some double cover $M$ of $N$, with no non-vanishing $D^*$- exact structure currents induces a non-empty contractible collection of non coorientable contact structures on $N$.
\end{thm}

\begin{rmk} Let us recall that a locally conformal symplectic (lcs) structure on a manifold $P$ is a $2$-form $\omega$ on $P$ such that for each point $p \in P$ there exists a neighborhood $U$ of $p$ in $P$ and a positive smooth function $f$ defined on $U$ such that the form $f \omega|_{U}$ is symplectic. The local functions $f$ are such that the exact forms $d(ln f)$ agree on overlaps and define a global closed form $\theta$, uniquely determined by $\omega$ and called the \emph{Lee form} of $\omega$. The lcs form $\omega$ is closed for the Lichnerowicz differential $D_\theta \beta = d\beta + \theta \wedge \beta$ associated to $\theta$ and, conversely, any $D_\theta$-closed non degenerate $2$-form is a lcs structure with Lee form $\theta$. \par

Now a lcs form, being a non degenerate $2$-form, admits a compatible almost complex structure and therefore an ample cone of bivectors whose cone of structure currents does not contain non-trivial $D^*_\theta$-exact currents. Conversely, provided the latter space is closed in the space of currents, one deduces existence of a lcs structure whose Lee form is $\theta$ from that of an ample cone structure that does not admit $D_\theta^*$-exact structure cycles. \par

It remains to prove that the collection of $D^*_\theta$-exact currents is a closed subspace of $\D_*(P)$. The proof is quasi-identical to that of \lref{bclosed}. Indeed, it suffices to observe that the manifold $P$ admits a covering $\pi : Q\to P$ for which $\pi^*\theta = df$ for some smooth function $f$ on $Q$ and to replace, in the above-mentioned lemma, the manifold $S^1 \times M$ by $P$, the manifold $\R\times M$ by $Q$, the $1$-form $dt$ by $\theta$ and the group $\Z$ by the group of Deck transformation of the covering $\pi : Q \to P$.\par

Altogether this proves that existence of a lcs form is equivalent to existence of an ample cone structure whose cone structure does not intersect at least one of the various spaces of $D_\theta$-exact forms for $\theta$ running through the space of closed $1$-forms on $P$.
\end{rmk}

\bibliography{x} 
\bibliographystyle{alpha}

\end{document}